\documentclass[10pt]{article}
\font\smallit=cmti10

\usepackage{amssymb,latexsym,amsmath,epsfig,amsthm} 
\usepackage{hyperref}
\hypersetup{
    colorlinks=true,
    linkcolor=blue,
    urlcolor =blue, 
    citecolor = magenta}
\usepackage[noabbrev, capitalize, nameinlink]{cleveref}
\usepackage{stmaryrd}
\usepackage{enumitem}
\setlist[enumerate,1]{label=(\alph*), ref=(\alph*)}
\usepackage{longtable}

\makeatletter

\renewcommand\section{\@startsection {section}{1}{\z@}
{-30pt \@plus -1ex \@minus -.2ex}
{2.3ex \@plus.2ex}
{\normalfont\normalsize\bfseries\boldmath}}

\renewcommand\subsection{\@startsection{subsection}{2}{\z@}
{-3.25ex\@plus -1ex \@minus -.2ex}
{1.5ex \@plus .2ex}
{\normalfont\normalsize\bfseries\boldmath}}

\renewcommand{\@seccntformat}[1]{\csname the#1\endcsname. }

\makeatother

\newtheorem{theorem}{Theorem}
\newtheorem{lemma}[theorem]{Lemma}
\newtheorem{proposition}[theorem]{Proposition}
\newtheorem{corollary}[theorem]{Corollary}

\newtheorem{thma}{Theorem}
\newtheorem{propa}[thma]{Proposition}
 
\theoremstyle{definition}
\newtheorem*{questionsx}{Questions}
\newenvironment{questions}
  {\pushQED{\qed}\questionsx}
  {\popQED\endquestionsx}
\newtheorem*{notatex}{Notation}
\newenvironment{notate}
  {\pushQED{\qed}\notatex}
  {\popQED\endnotatex}
\newtheorem*{defnx}{Definition}
\newenvironment{defn}
  {\pushQED{\qed}\defnx}
  {\popQED\enddefnx}
\newtheorem{exampx}[theorem]{Example}
\newenvironment{examp}
  {\pushQED{\qed}\exampx}
  {\popQED\endexampx}
\newtheorem{remx}[theorem]{Remark}
\newenvironment{rem}
  {\pushQED{\qed}\remx}
  {\popQED\endremx}

\newcommand{\ZZ}{\mathbb{Z}}
\newcommand{\FF}{\mathbb{F}}
\newcommand{\tr}{\operatorname{tr}}
\newcommand{\QQ}{\mathbb{Q}}
\renewcommand{\ss}{{\rm ss}}
\newcommand{\cmod}[1]{\ \mathrm{mod}\ #1}
\newcommand{\barr}[1]{\mkern 2mu\overline{\mkern-2mu#1\mkern-2mu}\mkern 2mu}
\renewcommand{\aa}{\mathfrak a}
\newcommand{\p}{{\mathsf p}} 
\newcommand{\g}{{\mathsf g}} 
\newcommand{\e}{{\mathsf e}} 
\newcommand{\into}{\hookrightarrow}
\newcommand{\lb}{\llbracket}
\newcommand{\rb}{\rrbracket}
\newcommand{\bb}{{\mathfrak b}}
\newcommand{\OO}{{\mathcal O}}
\newcommand{\mm}{{\mathfrak m}}
\newcommand{\PPP}{{\mathcal P}}
\newcommand{\NN}{\mathbb{N}}
\newcommand{\onto}{\twoheadrightarrow}
\newcommand{\tensor}{\otimes}
\newcommand{\E}{{\mathsf E}}
\renewcommand{\P}{{\mathsf P}}
\DeclareMathOperator{\rank}{{\mathrm rank}}
\newcommand{\mobius}{M\"obius}
\newcommand{\vdashp}{\vdash^{(p)}}
\newcommand{\eps}{\varepsilon}
\DeclareMathOperator{\charpoly}{CharPoly}
\DeclareMathOperator{\cchar}{{\rm char}}

\renewcommand{\geq}{\geqslant}

\renewcommand{\leq}{\leqslant}

\begin{document}

\begin{center}
\uppercase{\bf Deep congruences + the Brauer-Nesbitt theorem}
\vskip 20pt
{\bf Samuele Anni\footnote{The research of Samuele Anni is partially funded by the Melodia ANR-20-CE40-0013 project.}}\\
{\smallit Institut de Mathématiques de Marseille, Aix-Marseille Université, France}\\
{\tt samuele.anni@univ-amu.fr}\\
\vskip 10pt
{\bf Alexandru Ghitza}\\
{\smallit School of Mathematics and Statistics, University of Melbourne, Australia}\\
{\tt aghitza@alum.mit.edu}\\
\vskip 10pt
{\bf Anna Medvedovsky\footnote{Anna Medvedovsky was partially supported by NSF postdoctoral research fellowship DMS-1703834.}}\\
{\smallit Max-Planck-Institut für Mathematik, Bonn, Germany}\\
{\tt medvedov@post.harvard.edu}\\
\end{center}
\vskip 20pt

\centerline{\phantom{\smallit Received: , Revised: , Accepted: , Published: }}
\vskip 30pt

\centerline{\bf Abstract}

\noindent
We prove that mod-$p$ congruences between polynomials in $\ZZ_p[X]$ are equivalent to deeper $p$-power congruences between power-sum functions of their roots. This result generalizes to torsion-free $\ZZ_{(p)}$-algebras modulo divided-power ideals. Our approach is combinatorial: we introduce a $p$-equivalence relation on partitions, and use it to prove that certain linear combinations of power-sum functions are $p$-integral.  We also include a second proof, short and algebraic, suggested by an anonymous referee.
As a corollary we obtain a refinement of the Brauer-Nesbitt theorem for a single linear operator, motivated by the study of Hecke modules of mod-$p$ modular forms.

\pagestyle{myheadings}
\markright{\smallit Deep congruences\hfill}
\thispagestyle{empty}
\baselineskip=12.875pt
\vskip 30pt


\section{Introduction}

\subsection{The basic module-theoretic question} {Let $p$ be a prime.}
For a finite free $\ZZ_p$-module $M$ with an action of a linear operator~$T$, how much information does one need to know about the traces of $\ZZ_p[T]$ acting on~$M$ to know the structure of the semisimplification of $M \otimes \FF_p$ as an $\FF_p[T]$-module? 

Certainly knowing $\tr (T^n | M)$ as an element of $\ZZ_p$ for enough $n$ is plenty: the Brauer-Nesbitt theorem --- or in this one-parameter case, even simply linear independence of characters 
(see \hyperref[brauernesbittsec]{Appendix})
 --- tell us that these traces determine $(M \otimes \QQ_p)^\ss$, so that they determine the multiset of eigenvalues of $T$ on $M$ in characteristic zero, and hence in characteristic $p$.
 But this very precise characteristic-zero 
 information
 is more than we need: we merely want to understand $M$ modulo~$p$. 

On the other hand, knowing all the $\tr (T^n | M)$ modulo $p$ is not enough to determine $M \otimes \FF_p$. Indeed, if $M$ has rank $p$ and $T$ acts on $M$ as multiplication by a scalar $\alpha$ in $\ZZ_p$ then for every $n \geq 0$ we have \mbox{$\tr (T^n | M) = p \alpha^n \equiv 0 \cmod{p}$}, and we cannot recover $\alpha \cmod{p}$ from this trace data.

Since knowing $\tr (T^n | M)$ in $\ZZ_p$ is too much and knowing $\tr (T^n | M)$ modulo $p$ is not enough, one can ask for some kind of in-between criterion depending on $\tr (T^n|M)$ modulo \emph{powers} of $p$. This is the purpose of the present text: we precisely describe the exact depth of the $p$-adic congruence that the $\tr (T^n |M)$ must satisfy in order to pin down $M \otimes \FF_p$ up to semisimplification, and nothing more. In particular, we prove the following theorem.

\begin{thma}[{see \cref{moduletheorem}}] \label{theorema}
Let $M$ and $N$ be two finite free $\ZZ_p$-modules of the same rank $d$, each with an action of an operator $T$. Then
$\barr M^\ss \simeq \barr N^\ss$
as modules over $\FF_p[T]$ if and only if for every $n$ with $1 \leq n \leq d$
we have
\begin{equation*}
  \tr (T^n | M) \equiv  \tr (T^n | N) \mod{pn}.
\end{equation*}
\end{thma}

Here $\barr M$ and $\barr N$ are the $\FF_p[T]$-modules $M \otimes \FF_p$ and $N \otimes \FF_p$, respectively, and $\barr M^\ss$ and $\barr N^\ss$ refer to their semisimplifications. We highlight a few observations. 
\begin{itemize}
\item Since every prime except $p$ is a $\ZZ_p$-unit, congruence modulo $pn$ is the same as congruence modulo $p^{1 + v_p(n)}$, where $v_p: \QQ_p \to \ZZ$ is the $p$-adic valuation, normalized so that $v_p(p) = 1$.
\item \cref{theorema} completely resolves our example with $T = \alpha$ acting on $M = \ZZ_p^{\oplus p}$: knowing $\tr (T^p | M) = p \alpha^p$ \mbox{modulo~$p^2$} 
is tantamount to knowing $\alpha^p$ modulo $p$, which in turn determines $\alpha$ modulo $p$ uniquely. Yet this information is not enough to pin down $\alpha$ in $\ZZ_p$. 

\item The ``only if" direction of \cref{theorema} is trivial when all the eigenvalues of $M$ and $N$ are in $\ZZ_p$. Indeed, $\barr M^\ss \simeq \barr N^\ss$ implies that  eigenvalues of $M$ and $N$ pair by  mod-$p$ congruence.  But the $(p^{k})^{\rm th}$ powers of two mod-$p$-congruent elements  of $\ZZ_p$ are congruent modulo $p^{k + 1}$ (see~\cref{vpbase}); the deeper congruence claim follows. 
Thus the heart of \cref{theorema} is the ``if" direction. 

\item \cref{theorema} generalizes to valuation rings of $p$-adic fields that are not too ramified: see \cref{moduletheorem}. 
\end{itemize} 

The proof of \cref{theorema}, combinatorial in nature, follows from the slightly more general \cref{theoremb}, described in the next subsection.

\paragraph*{NB.} \label{NB}
An anonymous referee of this document
suggested a much simpler proof of \cref{theorema} than the one we present; see \cref{refproof}.
We still believe that our notion of $p$-equivalence for partitions --- and in particular \cref{propc} (the proof of which given here is due to Ira Gessel) --- 
used in the proof of \cref{theoremb}%
, as well as the observation in \cref{dpatp} (which we have not seen in the literature), have something to offer%
, so we present them here. It is also possible to prove \cref{theorema} purely algebraically, drawing inspiration from the proof of the characteristic-$p$ refinement of the trace version of Brauer-Nesbitt theorem (see \cref{brauernesbitt}\ref{tracecomp}) plus some algebra. The dedicated reader may find this third proof in our \href{https://arxiv.org/abs/2207.07108v1}{first Arxiv~draft}.

\subsection{The combinatorial perspective} Viewing \cref{theorema} as a combinatorial statement about deep congruences between power-sum symmetric functions 
implying simple congruences between corresponding elementary symmetric functions permits more generality. Let $A$ be a torsion-free $\ZZ_{(p)}$-algebra; for the purposes of this introduction only, we also assume that $A$ is a domain. Let $\aa \subset A$ be a \emph{divided-power ideal} --- see \cref{dpideal} for details and discussion, but in short, we must have $a^p \in p \aa$ for any $a \in \aa$. For a monic polynomial $P \in A[X]$, write $\barr P$ for the image of $P$ in $(A/\aa)[X]$ and $\p_n(P)$ for the $n^{\rm th}$ power-sum symmetric function of the roots of $P$ --- see Notation in \cref{lambda} for more and for the non-domain case. The following combinatorial theorem is a generalization of \cref{theorema}. 
\begin{thma}[{see \cref{fullalgebrathm}}] \label{theoremb}
  Let $P, Q$ be monic polynomials in $A[X]$. Then
  \begin{align*}
    \barr P = \barr Q\text{ in }(A/\aa)[X]
    \iff & \p_n(P) \equiv \p_n(Q)\text{ modulo }n \aa \\
         & \text{for }1 \leq n \leq \max\{ \deg P , \deg Q\}.
  \end{align*}
\end{thma}

In particular, here we do not require $P$ and $Q$ to be of the same degree; nor do we require $\aa$ to be prime (nor indeed $A$ to be a domain). 

The proof of \cref{theoremb} uses combinatorial theory of symmetric functions, specifically, formulas that express elementary symmetric functions in terms of power-sum functions and vice versa. Both directions of these formulas are sums indexed by partitions; for the ``if" direction, we introduce a new equivalence relation called \emph{$p$-equivalence} on the space of partitions to break up the sum: see \cref{pequivalencesec} for exact definitions --- but, for example, partitions $(6, 2)$, $(3, 3, 2)$, $(6, 1, 1)$, and $(3, 3, 1, 1)$ are all $2$-equivalent. 
The raison d'\^etre of $p$-equivalence is the following proposition. 

\begin{propa}[{see \cref{gesselprop}}]\label{propc}\ \\
Fix a partition $\lambda$ of an integer $n$. Write $C_\lambda$ for the set of partitions of $n$ that are $p$-equivalent to $\lambda$. Then the symmetric function 
\begin{equation*}
  \g_\lambda := \sum_{\mu \in C_\lambda} \frac{(-1)^\mu}{z_\mu} \p_\mu
  \quad\text{has coefficients in $\ZZ_{(p)}$.}
\end{equation*}
\end{propa} 

Here $(-1)^\mu$ is the sign in $S_n$ of any permutation $\sigma$ with cycle structure $\mu$, and $n!/z_\mu$ is the size of the $S_n$-conjugacy class of such a $\sigma$ (\cref{partdef}); the symmetric function $\p_\mu$ is the product of power-sum functions associated to the parts of $\mu$ (\cref{lambda}). For context, the elementary symmetric function $\e_n$ is the sum of the $\g_\lambda$ as $\lambda$ runs through a set of representatives of the $p$-equivalence classes (see \cref{glamdef} for details).

The elegant proof of \cref{propc} that we present  in \cref{gesselsec}, which relies on the $p$-integrality of the Artin-Hasse series, is due to Ira Gessel. We hope that the $p$-equivalence relation may be of independent interest in the study of partitions. 
 
\subsection{A generalization to virtual modules} The final result that we highlight in this introduction is a corollary of \cref{theorema}. 

\begin{corollary}
  \label{zpquotisom}
  Let $M_1,M_2,N_1,N_2$ be free $\ZZ_p$-modules of finite rank, each with\linebreak an action of an operator $T$.
  Suppose we have fixed $T$-equivariant embeddings\linebreak
  $\iota_1: \barr{N_1}\into\barr{M_1}$ and $\iota_2: \barr{N_2}\into\barr{M_2}$ and consider the quotients
  \begin{equation*}
    W_1:=\barr{M_1}/\iota_1(\barr{N_1}),\qquad W_2:=\barr{M_2}/\iota_2(\barr{N_2}).
  \end{equation*}
  Then $W_1^\ss \simeq W_2^\ss$ as $\FF_p[T]$-modules if and only if for every $n\geq 0$ we have
  \begin{equation*}
    v_p\big(\tr (T^n | M_1) -  \tr (T^n | N_1) - \tr (T^n | M_2) + \tr(T^n | N_2)\big) \geq 1 + v_p(n).
\end{equation*}
\end{corollary}

The essential point is that we do not assume that there are embeddings $N_i\into M_i$ over $\ZZ_p$, but only after base change to $\FF_p$. \cref{zpquotisom} is the form of the result that we use in \cite{refdim} to study the Hecke module structure on certain quotients of spaces of mod-$p$ modular forms. This is the motivating application of the present work, which we describe briefly below. 

\subsection{Motivating application to modular forms}

For $N$ prime to $p$ and $k \geq 2$, write $M_k(Np, \ZZ_p)$ for the space of classical modular forms of weight $k$ and level $Np$, viewed via the $q$-expansion map as a finite rank free $\ZZ_p$-submodule of $\ZZ_p\lb q \rb$. Let $M_k(Np, \FF_p)$ denote the image of $M_k(Np, \ZZ_p)$ in $\FF_p\lb q \rb$.  For $k \geq 4$, multiplication by the Eisenstein series $E_{p-1}$ normalized to be in $1 + p\ZZ_p\lb q \rb$  induces an embedding  
$ M_{k- p +1}(Np, \FF_p) \into M_{k}(Np, \FF_p);$
let $$\mbox{$W_k(Np) := M_{k}(Np, \FF_p)/M_{k-p+1}(Np, \FF_p)$}$$ denote the quotient.
 In \cite{refdim} we use \cref{zpquotisom} to prove that, for $p \geq 5$,
\begin{equation}\label{bp} W_k(Np)^\ss[1] \simeq W_{k + 2}(Np)^\ss
\end{equation}
as modules for the Hecke algebra generated by the action of Hecke operators $T_m$ for $m$ prime to~$Np$ (this is the \emph{anemic} or \emph{shallow} Hecke algebra).
The notation $W[1]$ stands for the Hecke module given by the vector space $W$ on which $T_m$ acts as $mT_m$ for all $m$ prime to~$Np$. We also refine \eqref{bp} to account for the action of the Atkin-Lehner involution at $p$ --- the main motivation for \cref{theorema}.

\paragraph*{Leitfaden.} \Cref{statementsec,combosec,onlyifsec,ifsec} are devoted to the proof of \cref{theoremb}. In \cref{statementsec}, we state \cref{fullalgebrathm}, the most general version of \cref{theoremb}, after a detailed discussion of the divided-power property of an ideal. In \cref{combosec} we collect and at times slightly extend a number of well-known results about symmetric functions, $p$-valuations of multinomial coefficients, and the $p$-integrality of the Artin-Hasse exponential series. We include proofs, both for completeness and because we hope that the motivating application will lure readers less familiar with combinatorics. In \cref{onlyifsec,ifsec} we prove the two directions of \cref{fullalgebrathm}; in particular, \cref{ifsec} is the heart of our main work here. In \cref{secondproofsec}, we return to the module-theoretic \cref{theorema} and deduce it from \cref{fullalgebrathm}. In the same section we also prove \cref{zpquotisom}.

\section{Statement of the main theorem}\label{statementsec}

\subsection{A bit of symmetric function notation}
For any ring $B$ and monic polynomial $P \in B[X]$ of degree $d$, let $\e_n(P)$ be the $X^{d-n}$-coefficient of $P$ scaled by $(-1)^n$. If $B$ is a domain, then $P$ determines $d$ roots $\alpha_1, \ldots, \alpha_d$ in some integral extension of $B$, and  $\e_n(P)$ is the $n^{\rm th}$ elementary symmetric function in the $\alpha_i$: namely, $$\e_n(P) = \displaystyle\sum_{ 1 \leq i_1 < i_2 < \cdots < i_n \leq d} \alpha_{i_1} \cdots \alpha_{i_n}.$$ Write $\p_n(P):=\sum_{i=1}^d \alpha_i^n$ for the $n^{\rm th}$ power-sum function of the roots of $P$. For a general $B$, Newton's identities \cite[I.2.{$11^\prime$}]{macdonald} express $\p_n$ as an integer polynomial in $\e_1, \ldots, \e_d$, thus defining $\p_n(P)$, or see \cref{lambda} below.

\subsection{Divided-power ideals in torsion-free \texorpdfstring{$\ZZ_{(p)}$}{Zp}-algebras}\label{dpideal}
Fix a torsion-free $\ZZ_{(p)}$-algebra\footnote{Recall that $\ZZ_{(p)}\subseteq\QQ$ is the subring of rationals that can be expressed as $\frac{a}{b}$ where $p\nmid b$.} $A$; in particular, $A$ embeds into $A[\frac{1}{p}] = A \otimes_{\ZZ_{(p)}} \QQ$. We say that an ideal $\aa$ of $A$ \emph{satisfies the divided-power property at} some $k \geq 1$ if $a \in \aa$ implies that  $\displaystyle {a^k}/{k!}$ is also in $\aa$. Since $A$ is $\ZZ$-torsion free and a $\ZZ_{(p)}$-algebra, this last condition may be reformulated: indeed, we have  
$$\mbox{$\displaystyle \frac{a^k}{k!}$ is in $\aa$} \iff \mbox{$a^k$ is in $k! \aa$} \iff \mbox{$a^k$ is in $p^{v_p(k!)} \aa$}.$$
An ideal $\aa$ that satisfies the divided power property for all $k \geq 1$ will be called a \emph{divided-power ideal}.
This concept plays a key role in the theory of crystalline cohomology, where $\aa$ satisfying the above condition exactly means that the maps $\gamma_k\colon\aa\to A$ given by $\gamma_k(a)=\frac{a^k}{k!}$ define a \emph{divided-power structure} on $\aa$~\cite[\S{}3]{berthelotogus}.

In a torsion-free $\ZZ_{(p)}$-algebra, satisfying the divided-power property at $p$ only is equivalent to being a divided-power ideal, as the following proposition shows.
\begin{proposition} \label{dpatp}
For an ideal $\bb$ in a commutative ring $B$, the following are equivalent
\begin{enumerate}
\item\label{atalln} For all $n \in \ZZ^+$ and all $a \in \bb$, we have $a^n \in p^{v_p(n!)} \bb$.
\item\label{atp} For all $a \in \bb$ we have $a^p \in p \bb$.
\end{enumerate}
\end{proposition}

\begin{proof}
The implication \ref{atalln} $\implies$ \ref{atp} is immediate given that $v_p(p!) = 1$. Suppose now that \ref{atp} is satisfied. First we show that \ref{atalln} is true for $n = p^k$ by induction on~$k$. The case $k = 0$ is trivial and $k = 1$ is exactly \ref{atp}. Suppose now \ref{atalln} is true for $n = p^k$ for some $k \geq 1$. Note that $$v_p(p^{k + 1}!) = p^k + p^{k-1} + \cdots + 1 = p v_p(p^k!) + 1.$$
For any $a \in \bb$, there exists a $b \in \bb$ so that $a^{p^k} = p^{v_p(p^k!)} b.$ Therefore
$$a^{p^{k + 1}} = (a^{p^k})^p = \big(p^{v_p(p^k !)} b\big)^p = p^{p v_p(p^k !)} b^p.$$
Since $b \in \bb$, by the \ref{atp} assumption we have $b^p \in p \bb$. Therefore
$$a^{p^{k +1}} \in p^{p v_p(p^k !) + 1} \bb = p^{v_p(p^{k + 1}!)} \bb,$$
as desired. 

Now for general $n \geq 1$, write $n$ in base $p$ as $n = n_k p^k + \cdots + n_1 p + n_0$, with \mbox{$n_i \in \{0, \ldots, p - 1\}$} for $i = 0, \ldots, k$. Fix $a \in \bb$ again. Since we have shown that for every $i$ we have
$a^{p^i} \in p^{v_p(p^i !)} \bb,$
we have
$a^{n_i p^i} \in p^{n_i v_p(p^i !)} \bb,$
so that
$a^n \in p^{ \sum_{i = 0}^k n_i v_p(p^i !)} \bb.$
The desired statement follows by observing that
$$\sum_{i = 0}^k n_i v_p(p^i !) = \sum_{i = 0}^k n_i \frac{p^{i} - 1}{p - 1} = \frac{n - \sum_{i = 0}^k n_i}{p-1} = v_p(n!),$$
where the last equality follows from a refinement of Legendre's formula on valuations of $n!$ (for a convenient exposition of this refinement, see \cite{romagny}).
\end{proof}

\begin{corollary} \label{dpredef}\ \\
The ideal $\aa \subset A$  is a divided-power ideal if and only if $a^p \in p\aa$ for every $a \in \aa$. 
\end{corollary} 

In fact, it suffices to check the condition of \cref{dpredef} on generators.

\begin{proposition} \label{dpgen} Let $S \subseteq A$ be a subset. Then the ideal $\aa$ generated by $S$ is a divided-power ideal if and only if ${a^p} \in p \aa$ for every $a \in S$. 
\end{proposition} 

\begin{proof} 
It suffices to show that for $a_1, a_2$ in $S$, and $b_1, b_2$ in $A$, if $a_1^p$ and $a_2^p$ are both in $p \aa$, then so is $(b_1 a_1 + b_2 a_2)^p$. We expand 
$$(b_1 a_1 + b_2 a_2)^p = b_1^p a_1^p + \sum_{k = 1}^{p-1} \binom{p}{k} b_1^k a_1^k b_2^{p-k} a_2^{p-k} + b_2^p a_2^p.$$
The first and last terms are in $p\aa$ by assumption; the middle terms  because \mbox{$p \mid \binom{p}{k}$.}
\end{proof}

\begin{corollary}\label{dpinduct}
If $\aa \subset A$ is a divided-power ideal, then so is $\aa\bb$ for any ideal $\bb \subseteq A$.
\end{corollary}

\begin{proof} For $a \in \aa$, $b \in \bb$ we have $(ab)^p = a^p b^p \in (p\aa) b^p \subseteq p(\aa \bb)$. Then  \cref{dpgen}. 
\end{proof}

\subsection{Divided-power ideals in \texorpdfstring{$p$}{p}-adic DVRs}
Recall that $v_p: \QQ_p \to \ZZ$ denotes the usual $p$-adic valuation, normalized so that $v_p(p) =1$. Let  $\OO$ be the ring of integers in a finite extension of $\QQ_p$, so that $v_p$ extends uniquely to $\OO$. Then $\OO$ is a torsion-free $\ZZ_{(p)}$-algebra and a complete DVR, so we will refer to such an $\OO$ as a \emph{$p$-adic DVR}. Any results for $p$-adic DVRs below also hold for localizations of rings of integers of number fields at prime ideals above~$p$; these are local torsion-free $\ZZ_{(p)}$-algebras whose completions are $p$-adic DVRs in the sense above, 
with completion establishing a one-to-one correspondence of ideals preserving the divided-power property.

\begin{lemma}  \label{OO} \ \\
 An ideal $\aa$ of a $p$-adic DVR is a divided-power ideal if and only if 
$v_p(\aa) \geq \frac{1}{p-1}$.
\end{lemma}

\begin{proof}
Let $a \in \aa$ be a generator, so that $v_p(a) = v_p(\aa)$. By \cref{dpgen}, the ideal $\aa$ is a divided-power ideal if and only if $a^p \in p \aa$, which happens in our $p$-adic DVR setting if and only if $$p v_p(a) = v_p(a^p) \geq v_p(p \aa) = 1 + v_p(a);$$ in other words, if and only if $v_p(a) \geq \frac{1}{p-1}$. 
\end{proof}

\begin{corollary} \label{Oe} Let $\mm$ be the maximal ideal of a $p$-adic DVR $\OO$. Let $e$ be the ramification degree of~$\mm$ over $p$.
Then $\mm$ is a divided-power ideal of $\OO$ if and only if $e \leq p -1$. In particular, $(p)$ is a divided-power ideal of $\ZZ_p$. 
\end{corollary}
\begin{proof}
Immediate from \cref{OO} as $v_p(\mm) = \frac{1}{e}$ in this setting. 
\end{proof}

\subsection{Statement of the main theorem}
We are ready to state the fullest version of \cref{theoremb}.

\begin{theorem} \label{fullalgebrathm}\ Let $A$ be a torsion-free $\ZZ_{(p)}$-algebra and $\aa$ a divided-power ideal, and let $P, Q$ be monic polynomials in $A[X]$. Then the following are equivalent:
\begin{enumerate}[itemsep = 0pt, topsep = 0pt]
\item\label{11}  $\e_n(P) \equiv \e_n(Q) \cmod{\aa}$ for every $n \geq 1$;
\item \label{22} $\e_n(P) \equiv \e_n(Q) \cmod{\aa}$ for every $n$ with $1 \leq n \leq \max\{ \deg P, \deg Q\}$;
\item \label{33} $\p_n(P) \equiv \p_n(Q) \cmod{n \aa}$ for every $n \geq 1$;
\item \label{44} $\p_n(P) \equiv \p_n(Q) \cmod{n \aa}$ for every $n$ with $1 \leq n \leq \max\{ \deg P, \deg Q\}$.
\end{enumerate}
\end{theorem}

\begin{rem} \label{fullalgthmext}
We do not require $\deg P = \deg Q$ here. In fact, since the statement  $\deg P = \deg Q$ is the same as the congruence
    $\p_n(P) \equiv \p_n(Q) \mod{n \aa}\quad\text{for }n = 0$,
  we may if we like  replace $n \geq 1$ with $n \geq 0$ in \ref{33} and~\ref{44} at the price of adding the condition $\deg P = \deg Q$ in \ref{11} and \ref{22}. In this case, we may add a fifth equivalent statement to \cref{fullalgebrathm}: 
\begin{enumerate}
\addtocounter{enumi}{4}
\item \label{55} $\barr P = \barr Q$ in $(A/\aa)[X]$.\qedhere
\end{enumerate}
\end{rem}

\begin{examp}\label{kickass}
Let $p = 2$ and $A = \ZZ_p$; given the polynomials $P = X^2 + X + 3$ and $Q = X^4 + 3X^3 + 5X^2 + 2X + 6$, we have
\begin{center}
\begin{tabular}{r||r|r|r|r||c|c@{\qquad\;\;\;}}
$n$ &  $\e_n(P)$ & $\e_n(Q)$ & $\p_n(P)$ & $\p_n(Q)$ & $v_2\big( \p_n(Q) - \p_n(P)\big)$ & $1+ v_2(n)$\\
\cline{1-7}
0  & 1 & 1 & 2 & 4 & 2 & $\infty$\\
1 & $-1$ & $-3$ & $-1$& $-3$ & 1 & 1
\\
2 & 3 & 5 & $-5$ &$-1$ & 2 & 2
\\
3 & 0 & $-2$ & 8& 12& 2 & 1
\\
4 & 0 & 6 &  7 &$-49$ &3 & 3
\\
5  & 0 & 0& -31 & 107 & 1 & 1
\\
6  & 0 & 0&10 &$-94$ & 3 & 2
\\
7& 0 & 0& 83& $-227$& 1& 1
\\
8 & 0 & 0&$-113$ &1231 &6 & 4
\\
9 & 0 & 0&$-136$ &$-3012$& 2 & 1
\\
10  & 0 & 0&475 &3899 &5 & 2
\\
11  & 0 & 0&$-67$ &2263& 1 & 1
\\
12 & 0 & 0&$-1358$ &$-27646$ & 4&3
\\
13& 0 & 0& 1559 & 81897 &1 & 1
\\
14  & 0 & 0&2515 &$-135381$& 3 & 2
\\
15  & 0 & 0&$-7192$& 38372& 2 & 1
\\
16  & 0 & 0&$-353$ &563871 &10 & 5
\end{tabular}
\end{center}

From matching up coefficients (or from the fact that $Q = (X^2 + 2X)P -(4X - 6)$), it is clear that $\e_n(P) \equiv \e_n(Q)$ modulo $2$ for every $n \geq 1$.
  In the table above, the last two columns illustrate \cref{fullalgebrathm}:
  \begin{equation*}
    v_2\big(\p_n(Q) - \p_n(P)\big) \geq 1 + v_2(n)\qquad\text{for }n\geq 1.\qedhere
  \end{equation*}
\end{examp}

We now give a skeleton proof of \cref{fullalgebrathm}. Technical details are postponed to \cref{ifsec,onlyifsec}.
\begin{proof}[{Proof of \cref{fullalgebrathm}}]
We clearly have \ref{33} $\implies$ \ref{44} and \ref{11} $\implies$ \ref{22}; moreover since $\e_n(P) = 0$ for $n > \deg P$ we have \ref{22}$ \implies$ \ref{11} as well, so that \ref{11} $\iff$ \ref{22}.

We  show that \ref{11} $\implies$ \ref{33} and  \ref{22}$\implies$\ref{44} by proving the following  in \cref{onlyifsec}.
\begin{proposition} \label{onlyifprop} Fix $N \geq 1$.\\ If $\e_n(P) \equiv \e_n(Q) \cmod{\aa}$ for all $1 \leq n \leq N$, then \mbox{$\p_N(P) \equiv \p_N(Q) \cmod{N\aa}$}.
\end{proposition}

We then show  \ref{33}$ \implies$ \ref{11} and \ref{44}$ \implies $ \ref{22} by proving the following in \cref{ifsec}.
\begin{proposition} \label{ifprop}
Fix $N \geq 1$. If $\p_n(P) \equiv \p_n(Q) \cmod{n\aa}$ for all $1 \leq n \leq N$, then \mbox{$\e_n(P) \equiv \e_n(Q) \cmod{\aa}$} for all $1 \leq n \leq N$.
\end{proposition}

Since we have shown that \ref{11} $\iff$ \ref{33} $\implies$ \ref{44} $\iff$ \ref{22} $\iff$ \ref{11}, we have a cycle and in particular deduce the equivalence of \ref{33} and \ref{44}.
\end{proof}

The divided-power property of the ideal $\aa$ is crucial to both directions of \cref{fullalgebrathm}. We illustrate this point by giving two counterexamples in the absence of this property. In both \cref{alexmagic1} and \cref{alexmagic2} below, let $\OO$ be the valuation ring of the field $\QQ_p(\alpha)$ where $\alpha = p^{\frac{1}{p}}$. Then the maximal ideal $\mm$ of $\OO$ is not a divided-power ideal (\cref{Oe}), having ramification degree $p$. In both cases, $P$ and $Q$ have the same degree $p$, so statements \ref{11} and \ref{22} of \cref{fullalgebrathm} are equivalent to the equality $\barr P = \barr Q$ in $\FF_p[X]$. 

\begin{examp}\label{alexmagic1}
Consider $P = X^p - \alpha X^{p-1}$ and $Q = X^p$. Then $P$ and $Q$, and hence their roots and their elementary symmetric functions are congruent modulo $\mm$. But $\p_p(P) = \alpha^p = p$ has $p$-valuation~$1$, and is not congruent to $\p_p(Q) = 0$ modulo $p \mm$, which has valuation $1 + \frac{1}{p}$.
Thus statements \ref{11} and \ref{22} of \cref{fullalgebrathm} hold but \ref{33} and \ref{44} do not.  
\end{examp} 

\begin{examp}\label{alexmagic2}
Consider $P = \big(X - (\alpha + p-1)\big)\big(X + 1\big)^{p-1}$ and $Q = X^p$. 
Then $P$ and $Q$ are \emph{not} congruent modulo~$\mm$: indeed, the roots of $P$ are units in $\OO$ whereas $Q$ has only zero as a root with multiplicity. 
 But we show that 
  \begin{equation*}
    \p_n(P) \equiv \p_n(Q) \equiv 0 \mod{n \mm}\quad\text{for }1 \leq n \leq p.
  \end{equation*}
  Indeed, for any $n \geq 1$, 
\begin{align} \begin{split} \label{alexmagic} 
\p_n(P) & = \big(\alpha + (p-1)\big)^n + (p-1)(-1)^n\\
&= \alpha^n + \sum_{i = 1}^{n - 1} \binom{n}{i} \alpha^i  (p - 1)^{n-i} + (p-1)^n + (p-1)(-1)^n\\
 &=(\mbox{terms divisible by $\alpha$}) + (p-1)^n + (p-1)(-1)^n.
\end{split} \end{align} 
Since $(p-1)^n + (p-1)(-1)^n \equiv (-1)^n - (-1)^n = 0$ modulo $p = \alpha^p$, we have $\p_n(P) \equiv 0$ modulo $\mm$.

If further $n = p$, then the summation term in \eqref{alexmagic} is divisible by $p\alpha = \alpha^{p+1}$, and the rest of the terms are $ \alpha^p + (p-1)^p + (p-1)(-1)^p$. If $p$ is odd, then
$$\alpha^p + (p-1)^p + (p-1)(-1)^p = p + (p-1)^p - (p-1) = (p-1)^p - (-1)  \equiv 0 \mod{p^2},$$
where the last congruence holds because $p-1 \equiv -1 \mod{p}$, so that their $p^{\rm th}$ powers are congruent modulo $p^2$ (see also \cref{vpbase} below). And if $p = 2$ then $$p + (p-1)^p + (p-1)(-1)^p = 2 + (-1)^2 + (1)(-1)^2 = 4.$$ In either case, 
$\p_p(P)$ is a sum of a term in $\mm^{p + 1}$ and a term in $\mm^{2p}$, so $\p_p(P) \in p\mm$, as required. 
Thus statement \ref{44} of \cref{fullalgebrathm} holds but \ref{11} and \ref{22} do not. One can show analogously that~\ref{33} also does not hold, as $v_p\big(\p_{2p}(P)\big)  = 1$.
\qedhere   
\end{examp}

\begin{questions}
  One can further ponder the relationship between the statements in \cref{fullalgebrathm}: 
\begin{itemize}
  \item Is there a direct proof of \ref{44} $\implies$ \ref{33}? The divided-power property or a similar assumption must play a role, as \cref{alexmagic2} above satisfies \ref{44} but not \ref{33}.
  \item Although \ref{44} does not imply \ref{11} or \ref{22} without the divided-power assumption (again, see \cref{alexmagic2} above), is it possible that \ref{33} does?\qedhere
\end{itemize}
\end{questions}

\medskip

\noindent The next three sections are devoted to the proof of \cref{fullalgebrathm}.

\section{Combinatorial preliminaries}\label{combosec}

\subsection{Partitions}\label{partdef}

A \emph{partition} $\lambda$ of an integer $n \geq 0$, denoted $\lambda \vdash n$, is a (finite or infinite) ordered tuple $(\lambda_1, \lambda_2, \ldots)$ with $\lambda_1 \geq \lambda_2 \geq \cdots \geq 0$ and $\sum_{i \geq 1} \lambda_i = n$. If the partition is infinite, only finitely many of the \emph{parts} $\lambda_i$ are nonzero. The number of nonzero parts of $\lambda$ is exactly the cardinality of $\{i \geq 1: \lambda_i > 0 \}$. There is a unique partition of $0$, namely $\varnothing \vdash 0$, the \emph{empty} partition. The following four definitions are standard.
\begin{itemize}
\item The \emph{weight} $| \lambda |$ of a partition $\lambda = (\lambda_1, \lambda_2, \ldots)$ is the number being partitioned: \mbox{$| \lambda | : = \sum_{i \geq 1} \lambda_i.$}
\item For $a \geq 1$, let $r_a(\lambda)$ be the number of times that $a$ appears as a part in $\lambda$.
\item For $\lambda \vdash n$, let $(-1)^{\lambda}$ be the sign of a permutation in $S_n$ with cycle structure $\lambda$. In other words, if $\lambda = (\lambda_1, \ldots, \lambda_k)$ with $\lambda_k > 0$, then $(-1)^{\lambda} = (-1)^{\sum_i (\lambda_i - 1)}.$
\item For $\lambda \vdash n$, let $z_\lambda := \prod_{a \geq 1} a^{r_a(\lambda)} r_a(\lambda)!$ be the order of the centralizer in $S_n$ of any permutation of cycle structure $\lambda$, so that
${n!}/{z_\lambda}$ is the number of permutations of $n$ with cycle structure~$\lambda$.  Accordingly, $z_\varnothing = 1$.
\end{itemize}

For $n \geq 0$, let $\PPP_n$ be the set of partitions of $n$, and let $\PPP := \bigcup_{n \geq 0} \PPP_n$ be the set of all partitions, graded by weight. We can multiply two partitions as follows: for $\lambda \vdash n$ and $\mu \vdash m$, let $\lambda \mu$ be the partition of $m + n$ whose parts are the union of the parts of $\lambda$ and $\mu$. This operation gives $\PPP$ the structure of a free abelian monoid on the set $\{(n): n \in \NN\}$ of partitions consisting of a single part. In particular, for any partition $\lambda \vdash n$ and any $k \geq 0$, we may consider the partition $\lambda^k \vdash kn$.

\begin{defn}
Let $p$ be a prime and $\lambda := (\lambda_1, \lambda_2, \ldots)$ a partition of $n \geq 0$. Define the \emph{$p$-valuation of~$\lambda$} by $v_p(\lambda): = \min_i \{ v_p(\lambda_i)\}.$ Note that $v_p(\lambda)$ is the greatest integer~$v$ with the property that we can express $\lambda$ as a $(p^v)^{\rm th}$ power: $\lambda = \mu^{p^v}$, where $\mu = (\lambda_1/p^v, \lambda_2/p^v, \ldots)$. Of course $v_p(\varnothing)= \infty$.
\end{defn}

\subsection{Ring of symmetric functions}\label{lambda}

Let $\Lambda_d$ be the ring of symmetric polynomials in $d$ variables $x_1, x_2, \ldots, x_d$ with integer coefficients: that is, $\Lambda_d$ consists of the $S_d$-invariants of $\ZZ[x_1, \ldots, x_d]$, where the symmetric group $S_d$ acts by permuting the variables. Then $\Lambda_d$ is a ring graded by degree: \mbox{$\Lambda_d = \bigoplus_{n \geq 0} \Lambda_d^n$}, where $\Lambda_d^n \subseteq \Lambda_d$ are the homogeneous symmetric polynomials in $x_1, \ldots, x_d$ of degree $n$. For any $d \geq d'$ we have a graded map $\Lambda_d \onto \Lambda_{d'}$ mapping $x_i$ to $x_i$ for $i \leq d'$ and sending $x_i$ with $i > d'$ to zero. This forms a compatible system of graded rings, and we take the so-called graded inverse limit to form the ring of symmetric functions: that is, $\Lambda^n:=\varprojlim_d \Lambda^n_d$ and $\Lambda := \bigoplus_{n \geq 0} \Lambda^n$. This somewhat fussy construction guarantees that every symmetric function in $\Lambda$ has finite degree.
For any ring $A$, let 
$\Lambda_A := \Lambda \tensor_\ZZ A$.

We now recall the definitions of some special symmetric functions and some general constructions.
\begin{itemize}
\item {\bf Elementary symmetric functions:} For $n \geq 0$, let $\e_{n, d} \in \Lambda_{d}^n$ be the \emph{$n^{\rm th}$ elementary symmetric polynomial}: $$\e_{n, d} := \sum_{ 1 \leq i_1 < i_2 < \cdots i_n \leq d} x_{i_1} \cdots x_{i_n},$$ and let $\e_n := \varprojlim_d \e_{n, d} \in \Lambda^n$ be the {$n^{\rm th}$ elementary symmetric function}. In particular $\e_0 = \e_{0, d} = 1$. One can check --- for example, see \cite[I.2.4]{macdonald} --- that 
\begin{equation}\label{pol}\Lambda = \ZZ[\e_1, \e_2, \ldots].\end{equation} 

\item {\bf Power-sum symmetric functions:} Similarly, for $n \geq 0$, let $\p_{n, d}\in\Lambda_d^n$ be the \emph{$n^{\rm th}$ power-sum polynomial}:
  \begin{equation*}
  \p_{n,d}:= \sum_{i = 1}^d x_i^n \in \Lambda_d^n.
  \end{equation*}
    For $n\geq 1$ we also let $\p_n: = \varprojlim_d \p_{n, d}\in \Lambda^n$ be the \emph{$n^{\rm th}$ power-sum function}. Note that $\p_{0, d} = d$, so that these do not interpolate and $\p_0$ is not defined.  One can check that $\Lambda_\QQ = \QQ[\p_1, \p_2, \ldots]$; see, for example, \cite[I.2.12]{macdonald}.

\item {\bf Symmetric functions depending on partition:}  We use the following standard notation: given a family of symmetric functions $\{f_n\}_{n \geq 1}$ --- for example, elementary or power-sum symmetric functions --- and a partition\linebreak
  $\lambda = (\lambda_1, \lambda_2, \ldots, \lambda_k)$, let $f_\lambda := f_{\lambda_1} f_{\lambda_2}\cdots f_{\lambda_k}$. In other words, we view $f$ as a map $(n) \mapsto f_n$ and extend it to a map of multiplicative monoids \mbox{$\PPP \to \Lambda$}. Note that $f_\varnothing = 1$. In particular, although $\p_0$ is undefined, we do have \mbox{$\p_{\varnothing} = \e_{\varnothing} = \e_0 = 1$}. We can also use the notation $f_\lambda$ for any tuple $\lambda$, not necessarily a partition. One can check that $\{e_\lambda\}_{\lambda \vdash n}$ is a $\ZZ$-basis for $\Lambda^n$ and $\{\p_\lambda\}_{\lambda \vdash n}$ is a $\QQ$-basis for $\Lambda^n_\QQ$.
\end{itemize}

Building on these, we introduce notation for a symmetric function evaluated at a polynomial. 
\begin{notate}
For a polynomial 
$Q = X^d + a_1 X^{d-1} + \cdots + a_d \in A[X]$ and $n \geq 0$, denote  by $$\e_n(Q) := \begin{cases} 1 & \mbox{if } n = 0\\
(-1)^n a_n  & \mbox{if }  1 \leq n \leq d \\
0 & \mbox{if } n > d.
\end{cases}$$
More generally, for any symmetric function $f$ and any monic polynomial $Q \in A[X]$, let $f(Q) \in A$ be defined as follows: first use \eqref{pol} to write $f$ as a polynomial in the~$\e_n$ and let $f(Q)$ be the result of plugging  $\e_n(Q)$ for $\e_n$ into that polynomial. If $A$ is a domain, this is equivalent to plugging in to $f$ the roots of $Q$ with multiplicity for the first $\deg Q$-many $x$s, and zeros for the rest. We extend this definition to $\p_0$, which is not a priori a symmetric function, by letting $\p_0(Q) :=\deg Q$. With this definition, the sequence $\{\p_n(Q)\}_{n \geq 0}$ satisfies an $A$-linear recurrence of order $\deg Q$, closely related to Newton's identities (see, for example, \cite[I.2.{$11^\prime$}]{macdonald}).
\end{notate}

\subsection{Combinatorial lemmas} Here we collect standard facts relating generating functions of various symmetric functions: see, for example, \cite[I.2]{macdonald}. For a set of positive integers $S \subseteq \NN$, let
\begin{equation*}
  \P_S(t) := \sum_{s \in S} (-1)^{s-1} \frac{\p_s}{s} t^s
\end{equation*}
be the weighted and signed power-sum generating function. 
Also set $\P(t): = \P_{\NN}(t)$.
On one hand, we can interpret the exponential of $\P_S(t)$ as a weighted sum of power-sum functions for partitions with parts restricted to $S$. The following proposition is standard for $S = \NN$;
 this formulation we learned from Gessel.

\begin{proposition}\label{iralog} Let $S \subseteq \NN$ be a set of positive integers.
  Then
  \begin{equation*}
  \exp \P_S(t) = \sum_{n = 0}^\infty \sum_{\substack{\lambda \vdash n \\{\rm parts\, in}\, S}} (-1)^\lambda \frac{\p_\lambda}{z_\lambda} t^n.
  \end{equation*}
\end{proposition}
\begin{proof}
\vspace{-15pt}
\begin{align*}
  \exp \P_S(t) &= \exp \left(\sum_{s \in S} (-1)^{s-1} \frac{\p_s}{s}t^s\right)\\  &= \prod_{s \in S} \exp \left( (-1)^{s - 1} \frac{\p_s}{s} t^s\right)\\
  & = \prod_{s \in S} \sum_{r_s = 0}^\infty \frac{1}{r_s!} (-1)^{r_s(s-1)} \frac{ \p_s^{r_s}}{s^{r_s}} t^{s r_s}\\  &= 
\sum_{(r_s) \in \NN^S} (-1)^{\sum_{s} r_s(s-1)} \frac{ \prod_{s} \p_s^{r_s}}{\prod_{s} r_s! s^{r_s}}t^{\sum_s s r_s} \\
& = \sum_{\lambda\ {\rm has\ parts\ in}\ S }   (-1)^\lambda \frac{\p_\lambda}{z_\lambda} t^{| \lambda |}.
\end{align*}
Here the sum in the penultimate line is over tuples of nonnegative integers $r_s$ indexed by elements of $S$ only finitely many of which are nonzero, and in the last line such a tuple is interpreted as a partition $\lambda$ all of whose parts are in $S$, with part $s$ appearing~$r_s$ times.
\end{proof}

On the other hand, for $S = \NN$ we can reinterpret $\exp \P_S(t)$ as the generating function for the elementary symmetric functions. Let $$\E(t) : = \sum_{k \geq 0} \e_k t^k = \prod_{i = 1}^\infty (1 + x_i t).$$ The remaining statements of this section are completely standard. 

\begin{proposition}\label{eexpp}
$\displaystyle \E(t) = \exp \P(t).$
\end{proposition}

\begin{proof}
We show that $\log \E(t) = \P(t)$:
\begin{align*}
\log \E(t) &:= \log \prod_{i = 1}^\infty (1 + x_i t) = \sum_{i = 1}^\infty \log (1 + x_i t)
 = \sum_{i = 1}^\infty \sum_{n = 1}^\infty (-1)^{n-1} \frac{(x_i t)^n}{n}\\
  &= \sum_{n = 1}^\infty (-1)^{n-1} \frac{t^n}{n} \sum_{i = 1}^\infty x_i^n
  = \sum_{n = 1}^\infty (-1)^{n-1} \frac{\p_n}{n} t^n = \P(t).\qedhere
\end{align*}
\end{proof}

\cref{eexpp} allows us to express $\e_n$ as a $\QQ$-linear combination of the $\p_\lambda$ for\linebreak 
$\lambda \vdash n$, and, conversely, $\p_n$ as a $\ZZ$-linear combination of $\e_\lambda$ over $\lambda \vdash n$: see \cref{enplam} and \cref{pnelam}.

\begin{corollary}[Expressing $\e_n$ in terms of $\p_\lambda$] \label{enplam} For all $n \geq 0$, we have
\begin{equation}\label{plamtoen}
\e_n =  \sum_{\lambda \vdash n} (-1)^\lambda \frac{\p_\lambda}{z_\lambda}. 
\end{equation}
\end{corollary}

For example, $\e_2 = \frac{\p_1^2 - \p_2}{2}$ and $\e_3 = \frac{\p_1^3 - 3\p_1 \p_2 + 2\p_3}{6}$. 

\begin{proof}
Combining \cref{eexpp} with \cref{iralog} for $S = \NN$ yields 
  \begin{equation*}
    \sum_{\lambda} (-1)^\lambda \frac{\p_\lambda}{z_\lambda} t^{|\lambda|} = \sum_{n = 0}^\infty \e_n t^n.
  \end{equation*}
The statement follows from considering the coefficients of $t^n$ on each side. \end{proof}

\begin{corollary}[Expressing $\p_n$ in terms of $\e_\lambda$]\label{pnelam} For $n \geq 1$, we have
\begin{equation}\label{entobnformula}
  \p_n  = (-1)^{n} \, n \sum_{\lambda \vdash n} \frac{ (-1)^{m}}{m} \binom{ m }{ r_1(\lambda), r_2(\lambda), \ldots} \e_\lambda,
\end{equation}
where $m:=r_1(\lambda)+r_2(\lambda)+\dots$ is the number of nonzero parts of the partition $\lambda$.
\end{corollary}

\begin{proof}
From \cref{eexpp} we have
\begin{align*}
  \sum_{n = 0}^\infty &(-1)^{n-1} \frac{\p_n}{n} t^n = \P(t) = \log \E(t) =  \log \left(1 + \sum_{k = 1}^\infty \e_k t^k\right)\\ 
  &= \sum_{m = 1}^\infty \frac{(-1)^{m-1}}{m} \left(\sum_{k = 1}^\infty \e_k t^k\right)^m
 = \sum_{m = 1}^\infty \frac{(-1)^{m-1}}{m} \sum_{1 \leq k_1, \ldots, k_m} \e_{k_1} \cdots \e_{k_m} t^{k_1 + \cdots + k_m},
\end{align*}
  where the last sum is over $m$-tuples $(k_1, \ldots, k_m)$ of positive integers. We can interpret such a tuple as a (badly ordered) partition $\lambda$ of $\sum k_i$ into $m$ parts, with $r_a(\lambda)$ of the $k_i$s equal to $a$ and $m = \sum_a r_a(\lambda)$. Moreover, each such partition $\lambda$ will arise from exactly $\displaystyle \binom{m}{r_1(\lambda), r_2(\lambda), \ldots }$ such $m$-tuples. Equating coefficients of $t^n$ on each side, we obtain, as desired, 
  \begin{equation*}
    \p_n = (-1)^{n-1} \, n \sum_{m \geq 1} \sum_{\substack{\lambda \vdash n \\ {\rm with}\ m\ {\rm parts}}} \frac{ (-1)^{m-1}}{m} \binom{m}{r_1(\lambda), r_2(\lambda), \ldots} \e_\lambda.\qedhere
\end{equation*}
\end{proof}

\subsection{A simple proof of \texorpdfstring{\cref{theorema}}{Theorem A}}\label{refproof}
An anonymous referee of this paper suggested a simpler proof of \cref{theorema}, which builds on the above discussion of the generating functions $\P$ and $\E$.

Let $M$, $N$ be free $\ZZ_p$-modules of rank $d$, each endowed with an action of an operator $T$.
 Writing $\p_n(M)$ and $\e_n(M)$ for the $n^{\rm th}$ power-sum and elementary symmetric function of the eigenvalues of $T$ on $M$, with the corresponding generating functions 
$$\P(M, t) := \sum_{n \geq 1} (-1)^{n-1} \frac{{\p_n(M)}}{n} t^n \in \QQ_p\lb t \rb \mbox{\quad and \quad} \E(M, t) : = \sum_{n \geq 0} \e_n(M) t^n \in \ZZ_p[t],$$ we note that we still have $\P(M, t) = \log \E(M, t)$ as in \cref{eexpp}, so that we may proceed as follows: 
\begin{align*} 
\barr M^\ss \simeq \barr N^\ss 
& \iff \mbox{for all $1\leq n \leq d$ we have}\ \e_n(M) \equiv \e_n(N) \mod{p} \\
&\iff \E(M, t) \equiv   \E(N, t) \mod{p\ZZ_p [t]} \\
&\iff \E(M, t) = \E(N, t) S(t) \mbox{ for some $S(t) \in 1 + t p\ZZ_p\lb t \rb$}\\
& \iff \log \E(M, t) = \log E(N, t) + \log S(t) \\
& \iff \P(M, t) = \P(N, t) + R(t) \mbox{ for some $R(t) \in tp \ZZ_p \lb t \rb$} \\
& \iff \mbox{for all $n \geq 1$ we have}\  \p_n(M) \equiv \p_n(N) \mod{n p}\\
& \iff \mbox{for all $n \geq 1$ we have}\ \tr(T^n | M) \equiv \tr(T^n | N) \mod{np}.
\end{align*}
Note that we used the fact that $\log$ maps $1 + t p \ZZ_p\lb t\rb$ onto $t p \ZZ_p\lb t \rb$.

The argument generalizes to the setting of \cref{theoremb}, with $\ZZ_p$ and $p$, respectively, replaced by torsion-free $\ZZ_{(p)}$-algebra $A$ and a divided-power ideal $\aa$ (see \cref{dpideal} for definitions), and the assumption $\rank M = \rank N$ relaxed. 

\subsection{\texorpdfstring{$p$}{p}-valuation lemmas}
Here we collect a few lemmas about $p$-valuations.  First, in light of the expression in \cref{pnelam} and our end goal, we need a formula for the $p$-valuation of multinomial coefficients. Let $r_1, \ldots, r_k$ be nonnegative integers, write $m = r_1 + \cdots r_k$, and let $p$ be any prime.  The following statement is due to Kummer for $k = 2$; see, for example \cite{romagny}. 
The generalization to any $k$ is immediate through the formula $$\binom{r_1 + \cdots + r_k }{ r_1, \ldots, r_k} = \binom{r_1 + \cdots + r_k }{ r_1}\binom{r_2 + \cdots + r_k }{ r_2} \cdots \binom{r_{k - 1} + r_k }{ r_{k -1}}$$ expressing multinomial coefficients in terms of binomial coefficients.

\begin{theorem}[Kummer, 1852]\label{kummercarry}\ \\ The $p$-valuation of the multinomial coefficient $\binom{m }{ r_1, \ldots, r_k}$ is the sum of the carry digits when the addition $r_1 + \cdots + r_k$ is performed in base $p$.
\end{theorem}

\begin{corollary}\label{kummercor}
For any $i$, we have
  $\displaystyle v_p(r_i) \geq v_p(m) -  v_p \left( \binom{m }{ r_1, \ldots r_k} \right).$
\end{corollary}

\begin{proof}
  Any end zero of $m$ base $p$ not corresponding to an end zero of $r_i$ base $p$ contributes to a carry digit of the base-$p$ computation $r_1 + \cdots + r_k =m$. Therefore, $v_p \left( \binom{m }{r_1, \ldots r_k} \right) \geq v_p(m) - v_p(r_i)$.
\end{proof}

The second statement we need (\cref{elamcongruence} below) is a partition version of the standard observation that the depth of the $p$-adic congruence of two integers increases upon taking $p^{\rm th}$ powers.

Recall that $A$ is a torsion-free $\ZZ_{(p)}$-algebra and $\aa \subset A$ is an ideal with a divided power structure.

\begin{lemma} \label{vpbase}
Suppose $x \equiv y \cmod{\aa}$ for some $x, y \in A$. Then 
\begin{enumerate}
\item for all $m \geq 0$ we have $x^{p^m} \equiv y^{p^m} \cmod{p^m \aa}$; more generally
\item for all $n \geq 0$ we have $x^n \equiv y^n \cmod{n \aa}.$
\end{enumerate}
\end{lemma}

\begin{proof}
Since $A$ is a $\ZZ_{(p)}$-algebra, the ideal $n \aa$ is the same as the ideal $p^{v_p(n)} \aa$. Thus it suffices to prove the first statement. 
  For $m = 1,$ write $y = x + a$ with $a \in \aa$. Then $$y^p - x^p = (x + a)^p - x^p = a^p + \sum_{i = 1}^{p-1} \binom{p }{ i} x^{p-i} a^i.$$  We show that each of the terms on the right-hand side is in $p\aa$. This is clear for each term in the summation because for $0 < i < p$ we have both $p \mid \binom{p }{ i}$ and $a^i \in \aa$. \cref{dpredef} tells us that $a^p\in p\aa$.
To prove the statement for $m > 1$ we proceed by induction using \cref{dpinduct}.
\end{proof}

\begin{corollary} \label{elamcongruence}
  Let $P, Q \in A[X]$ be polynomials, and $\{f_n\}_{n \geq 1}$ a family of symmetric functions. If $f_n(P) \equiv f_n(Q)$ modulo $\aa$ for all $n$, then for every partition $\lambda$
$$f_\lambda(P) \equiv f_\lambda(Q) \pmod{p^{v_p(\lambda)}\aa}.$$
\end{corollary}

\begin{proof}
Let $v = v_p(\lambda)$. By the definition of $p$-valuation of a partition (\cref{partdef}) there exists a partition $\mu$ so that $\lambda = \mu^{p^v}$. Therefore $$f_\lambda(P) = f_{\mu^{p^v}}(P) = f_\mu(P)^{p^v} \equiv_{p^v \aa} f_\mu(Q)^{p^v} = f_{\mu^{p^v}}(Q) = f_\lambda(Q),$$ where the middle congruence modulo $p^v \aa$ holds  by \cref{vpbase}.
\end{proof}

\subsection{Artin-Hasse exponential series} \label{artinhasse} We briefly recall the Artin-Hasse exponential series
\begin{equation*}
F(z)= \exp \left(\sum_{j = 0}^\infty \frac{z^{p^j}}{p^j} \right) =  1 + z + \frac{z^2}{2!} + \frac{z^3}{3!} + \cdots + \frac{z^{p-1}}{(p-1)!} + \frac{\big(\frac{(p-1)! + 1}{p}\big)z^p}{(p-1)!} + \cdots,
\end{equation*}
here viewed merely as a formal power series, a priori in $\QQ\lb z \rb$. In \cref{gesselsec} we will make use of the fact that $F(z)$ is actually $p$-integral (\cref{artinhassecor}); here we briefly review this well-known result. We follow the convenient expository notes \cite{lurie} of Jacob Lurie.
\begin{proposition} \label{mu}
We have $\displaystyle F(z) = \prod_{p \nmid d} \big(1 - z^d\big)^{-\frac{\mu(d)}{d}}.$
\end{proposition}

Here $\mu$ is the \mobius\ function, the multiplicative arithmetic function taking squarefree products of primes $p_1 \ldots p_k$ to $(-1)^k$ and other positive integers to $0$, and  satisfying the property 
\begin{equation}\label{mubasic} \sum_{d \mid n} \mu(d) = \begin{cases} 1 &  \mbox{if } n = 1\\ 0 & \mbox{otherwise.} \end{cases}\end{equation}

Before giving the proof of \cref{mu}, we need a lemma: 
\begin{lemma} \label{mup}
For prime $p$ we have $\displaystyle \sum_{d \mid n,\: p \nmid d} \mu(d) = \begin{cases} 1 & \mbox{if $n$ is a power of $p$} \\ 0 & \mbox{otherwise.} \end{cases}$
\end{lemma} 

\begin{proof} 
Quite generally if $f(n)$ is a multiplicative arithmetic function, then the function
$$\phi(n) := \sum_{d\mid n,\: p\nmid d} f(n)$$
is also multiplicative. Indeed, say  a divisor $d$ of $n$ is \emph{$p$-deprived} if $p\nmid n$. Then assuming \mbox{$\gcd(m,n) = 1$}, each $p$-deprived divisor of $mn$ is uniquely a product of a $p$-deprived divisor of $m$ and a $p$-deprived divisor of $n$, which are, in turn, relatively prime to each other. The fact that $f$ is multiplicative then allows the factorization $\phi(mn) = \phi(m) \phi(n)$. 

Now for the claim. Since $\mu$ is multiplicative, it suffices to check the claim for $n$ a power of $p$ and $n$ relatively prime to $p$. In the former case the claim is immediate; in the latter it follows from \eqref{mubasic}. 
\end{proof} 

\begin{proof}[{Proof of \cref{mu}}] 
We have 
\begin{align*}
\log \prod_{p \nmid d} \big(1 - z^d\big)^{-\frac{\mu(d)}{d}} \quad 
=\quad  &\sum_{p \nmid d}  \frac{\mu(d)}{d} \: \log \frac{1}{1 - z^d} 
 \quad=\quad 
 \sum_{p \nmid d}  \frac{\mu(d)}{d} \:\sum_{k \geq 1} \frac{z^{dk}}{k}\\
 =\quad &\sum_{n \geq 1} \frac{z^n}{n}\: 
 \sum_{d \mid n,\: p \nmid d} \mu(d) 
\,\!\:\:\: \quad=\quad 
  \sum_{ n = p^j,\: j \geq 0}
\frac{z^n}{n},
\end{align*} 
where the last equality follows from \cref{mup}. The claim follows. 
\end{proof} 

\begin{corollary} \label{artinhassecor} The Artin-Hasse exponential series $F(z)$ is in $\ZZ_{(p)}\lb z \rb$. 
\end{corollary}

\begin{proof} 
  The coefficients of $(1 - z^d)^{\pm 1/d}$ in the expression in \cref{mu} are algebraically generated by binomial coefficients $\binom{1/d}{ k }$, all in $\ZZ[\frac{1}{d}]$. Since all the $d$ are prime to $p$, the claim follows. 
\end{proof}

\section{Proof of \texorpdfstring{\cref{onlyifprop}}{Proposition 11}: \texorpdfstring{$\e_n$}{en} congruent implies \texorpdfstring{$\p_n$}{pn} deeply congruent}\label{onlyifsec}
Here we prove \cref{onlyifprop}. The proof uses the combinatorial expression from \cref{pnelam} for $\p_n$ in terms of $\e_\lambda$.

\begin{proof}[Proof of \cref{onlyifprop}] Let $P, Q \in A[X]$ be monic polynomials, fix $N \geq 1$, and suppose that $\e_n(P) \equiv \e_n(Q)$ modulo $\aa$ for all $n$ with $1 \leq n \leq N$. We seek to show that $\p_N(P) - \p_N(Q)$ is in $N \aa$.

From \cref{pnelam}
we have
  $$ \p_N(P) - \p_N(Q)  = (-1)^{N} \, N \sum_{m \geq 1} \sum_{\substack{\lambda \vdash N \\ {\rm with}\ m\ {\rm parts}}} \!\!\!\!\frac{ (-1)^{m}}{m} \binom{m }{ r_1(\lambda), r_2(\lambda), \ldots} \big(\e_\lambda(P) - \e_\lambda(Q)\big).$$
\cref{elamcongruence} for $f = \e$ tells us that our assumptions on the $\e_n$ imply that the difference $\e_\lambda(P) - \e_\lambda(Q) \in p^{v_p(\lambda)}\aa$ for each relevant $\lambda$. Therefore 
it suffices to show that for every $\lambda \vdash N$ with $m$ parts, 
  $$v_p(N) - v_p(m) +  v_p\left(  \binom{m }{ r_1(\lambda), r_2(\lambda), \ldots} \right)+ v_p(\lambda) \geq v_p(N),$$ or, equivalently, canceling $v_p(N)$ and using the definition of $v_p(\lambda)$, that for every $i$,
  $$-v_p(m) + v_p \left(\binom{m }{ r_1(\lambda), r_2(\lambda), \ldots} \right) + v_p\big(r_i(\lambda)\big) \geq 0.$$
But this is exactly \cref{kummercor}.
\end{proof}

Incidentally, although we know from the fact that the $\e_\lambda$ are a $\ZZ$-basis for $\Lambda$ in \eqref{entobnformula} that $\frac{n}{m}\binom{m }{ r_1, r_2, \ldots}$
is always integral --- here of course $r_1, r_2, \ldots $ is a sequence of nonnegative integers almost all zero, \mbox{$m = \sum r_i$} and $n = \sum i r_i$ --- it is not a priori obvious. But this integrality does follow from \cref{kummercor}. 

\section{Proof of \texorpdfstring{\cref{ifprop}}{Proposition 12}: \texorpdfstring{$\p_n$}{pn} deeply congruent implies \texorpdfstring{$\e_n$}{en} congruent}\label{ifsec}
Here we give the first, combinatorial, proof of the ``if" direction of \cref{fullalgebrathm}: we show that if the power sums of roots satisfy deep congruences, then elementary symmetric functions of the roots are (simply) congruent.

\subsection{\texorpdfstring{$p$}{p}-equivalent partitions}\label{pequivalencesec} We introduce an equivalence relation on the set $\PPP_n$ of partitions of an integer $n \geq 0$. 

\begin{defn} $\bullet$ If $\lambda$ and $\mu$ are in $\PPP_n$, we say that $\mu$ is a \emph{$p$-splitting} of $\lambda$ if $\lambda$ contains an instance of the part $pu$ for some $u \geq 1$, and $\mu$ is obtained from $\lambda$ by replacing~$pu$ with $p$ copies of part $u$. In other words, for every $u \in \NN$, the partition $(u)^p$ is a $p$-splitting of $(pu)$, and if $\mu$ is a $p$-splitting of $\lambda$, then $\mu \nu$ is a $p$-splitting of $\lambda \nu$.

$\bullet$ Let \emph{$p$-equivalence}, written $\sim_p$, be the equivalence relation generated by the $p$-splitting relation. For $\lambda \vdash n$ let  $C_{\lambda} = \{\mu \vdash n: \mu \sim_p \lambda\}$ denote the $p$-equivalence class of $\lambda$.

$\bullet$ Call a partition $\lambda$ of $n$ \emph{$p$-deprived} if none of its parts are divisible by $p$. The empty partition~$\varnothing$ is a $p$-deprived partition of $0$ for every $p$. Write $\lambda \vdash^{(p)} n$ for a $p$-deprived partition $\lambda$ of $n$. 
\end{defn}

  \begin{examp}\label{cur}
    Let $u\geq 1$ be prime to $p$ and let $r\geq 0$.
    Then the partition $(u)^r$ is $p$-deprived and 
    \begin{equation*}
      C_{(u)^r}=\big\{\lambda\vdash ur\colon \lambda\text{ has parts in }\{up^j\colon j\geq 0\}\big\}.\qedhere
    \end{equation*}
  \end{examp}

Every $p$-equivalence class has a unique $p$-deprived partition representative. We therefore have, for every $n \geq 0$, the following disjoint union: 
\begin{equation}\label{pdep} \PPP_n=\{\lambda \vdash n\} = \bigsqcup_{\lambda \vdash^{(p)} n} C_\lambda.
\end{equation}

\subsection{The contribution to \texorpdfstring{$\e_n$}{en} from a single \texorpdfstring{$p$}{p}-equivalence class}\label{glamdef}
Fix $n \geq 0$ and $\lambda \vdash n$. Let \begin{equation} \label{glam} \g_\lambda : = \sum_{\mu \sim_p \lambda} (-1)^\mu \frac{\p_\mu}{z_\mu},\end{equation}
so that in particular $\g_\varnothing = 1$. In other words, $\g_\lambda$ is the piece of the expression for~$\e_n$ from \eqref{plamtoen} that comes from all the partitions that are $p$-equivalent to $\lambda$. Because of~\eqref{pdep},
for any $n \geq 0$,
\begin{equation}\label{englam}
\e_n = \sum_{\lambda \vdashp n} \g_\lambda.
\end{equation}
To show that $\e_n(P) \equiv \e_n(Q) \cmod{\aa}$ in  \cref{ifprop}, it will therefore suffice to establish that \mbox{$\g_\lambda(P) \equiv \g_\lambda(Q) \cmod{\aa}$} for every $\lambda \vdashp n$.
But in fact we can break these up further:

\begin{lemma} \label{glammult}
Suppose $\lambda \vdash^{(p)} n$, $\mu \vdash^{(p)} m$ are partitions of $n, m \geq 0$ with no common parts. Then $$\g_{\lambda \mu} = \g_\lambda \g_\mu.$$
Thus for $\lambda \vdash^{(p)} n$,
$$\g_\lambda =   \prod_{u \geq 1,\ p \nmid u} \g_{(u)^{r_u(\lambda)}}.$$
\end{lemma}

\begin{proof}
Any two partitions $\lambda$ and $\mu$, whether disjoint or not, satisfy $\p_{\lambda \mu} = \p_\lambda \p_\mu$ and $(-1)^{\lambda  \mu} = (-1)^\lambda (-1)^{\mu}$. If $\lambda$ and $\mu$ have no parts in common, then $z_{\lambda \mu} = z_{\lambda} z_\mu$. And finally if both $\lambda$ and $\mu$ additionally have only prime-to-$p$ parts, then every $\nu \sim_p \lambda \mu$ factors uniquely as $\nu = \nu_{\lambda}  \nu_\mu$ with $\nu_\lambda \sim_p \lambda$ and $\nu_\mu \sim_p \mu$. The claim follows by the distributive property.
\end{proof}
Therefore rather than showing that $\g_{\lambda}(P) \equiv_\aa \g_{\lambda}(Q)$ for every $\lambda \vdash^{(p)} n$, it suffices to show that 
\begin{equation}\label{gur}\g_{(u)^r}(P) \equiv_\aa \g_{(u)^r}(Q)
\end{equation} for every $ur \leq n$ where $r \geq 0$ and $u \geq 1$ is prime to $p$.
We prove this in \cref{ifpropproof} after establishing a $p$-integrality result for the symmetric function $\g_\lambda$.

\subsection{\texorpdfstring{$p$}{p}-integrality of \texorpdfstring{$\g_\lambda$}{glambda}}\label{gesselsec}

First note that the signs $(-1)^\mu$ in the definition of $\g_\lambda$ are the same for every $\mu \sim_p \lambda$ for odd $p$. In other words,
\begin{lemma} \label{sign} If $p$ is odd, then  $\g_\lambda  = \displaystyle (-1)^\lambda \sum_{\mu \sim_p \lambda} \frac{\p_\mu}{z_\mu}$.
\end{lemma}

\begin{proof}
If $p$ is odd, then for any $u \geq 1$ and $j \geq 0$, the parity of $(u p^j)$ is the same as the parity of $(u)$ to the $p^j$ power: 
$$(-1)^{(u p^j)} = (-1)^{up^j-1} = (-1)^{u-1} = (-1)^{p^j(u-1)} = (-1)^{(u)^{p^j}}.$$
Then extend multiplicatively.
\end{proof}

From the definition in \eqref{glam} it is clear that $\g_\lambda$ is in $\Lambda_{\QQ}$. However, one can show that~$\g_\lambda$ is $p$-integral as a symmetric function. 

\begin{proposition} \label{gesselprop} For any partition $\lambda \vdash n \geq 0$, we have
$\g_\lambda$ in $\Lambda_{\ZZ_{(p)}}$.
\end{proposition}

The following elegant argument is due to Gessel. 

\begin{proof}
Since every equivalence class $C_\lambda$ has a unique representative with prime-to-$p$ parts, it suffices to consider $\g_\lambda$ for  $\lambda \vdash^{(p)} n$.
  By \cref{glammult}, it suffices to show that for any $u$ prime to $p$ and any $r \geq 0$, we have $\g_{(u)^r} \in \Lambda_{\ZZ_{(p)}}$. Equivalently, it suffices to show that for any $u$ prime to $p$, the generating function
\begin{equation}\label{gugessel} G_u(t) : = \sum_{r = 0}^\infty \g_{(u)^r}t^{ur}\end{equation}
for the sequence $\{g_{(u)^r}\}_{r \geq 0}$ is in $\Lambda_{\ZZ_{(p)}}\lb t \rb$. Recall that
  \begin{equation*}
  F(z)= \exp \left(\sum_{j = 0}^\infty \frac{z^{p^j}}{p^j} \right) \in \ZZ_{(p)}\lb z \rb
  \end{equation*}
  is the Artin-Hasse exponential series (\cref{artinhassecor}). 
 
 For $p$ odd, let $\eps_u = (-1)^{u - 1}$ be the sign of $(u{p^j})$ for $j \geq 0$ (\cref{sign}). Then
\begin{align}\begin{split}\label{gup} G_u(t)
& =
\exp \left( \sum_{j = 0}^\infty \frac{\eps_u \p_{u p^j} }{u p^j} t^{u p^j} \right)
  = \exp \left( \frac{\eps_u}{u} \sum_{j = 0}^\infty t^{u p^j}\frac{(x_1^{u p^j} + x_2^{u p^j} + \cdots)}{p^j} \right)\\
& = F(x_1^u t^u )^{\eps_u /u} F(x_2^u t^u)^{\eps_u /u} \cdots,
\end{split}\end{align}
  where  the first equality is \cref{iralog} for the set $S = \{up^j: j \geq 0\}$ (see \cref{cur}). 
  Since $F(x_i^u t^u)$ has coefficients in $\ZZ_{(p)}$
    and constant coefficient $1$, and since binomial coefficients $\binom{\eps_u/u }{ m}$ are in $\ZZ[\frac{1}{u}] \subset \ZZ_{(p)}$, each $F(x_i^u t^u)^{\eps_u /u}$ is in $\ZZ_{(p)}\lb x_i, t\rb$, so that $G_u(t)$ is in $\ZZ_{(p)}\lb t, x_1, x_2, \ldots \rb$. We already know it to be in $\Lambda_\QQ\lb t \rb$, so we conclude that $G_u(t) \in \Lambda_{\ZZ_{(p)}} \lb t \rb$, as desired.

It remains to consider $p = 2$. In this case, the sign of $(up^j)$ is $-1$ unless $j = 0$, in which case it is $1$ as $u$ is odd. Therefore, for $p = 2$,
\begin{align} 
  \nonumber G_u(t)
&= \exp \left(\frac{2 t^u \p_u}{u} -  \sum_{j = 0}^\infty \frac{t^{u p^j}}{u p^j} \p_{u p^j} \right)\\
  &= \left(\sum_{k = 0}^\infty \frac{2^k}{u^k k!}\, \p_u^k\, t^{uk} \right) \, F(x_1^u t^u )^{-1/u} F(x_2^u t^u)^{-1 /u} \cdots.\label{gu2}
\end{align}
To conclude that $G_u(t) \in \Lambda_{\ZZ_{(p)}}\lb t \rb$ for $p = 2$, we note that
\begin{equation*}
v_2(k!) = \left\lfloor \frac{k}{2} \right\rfloor + \left\lfloor \frac{k}{2^2} \right\rfloor + \cdots < \sum_{i = 1}^\infty \frac{k}{2^i} = k = v_2(2^k),
\end{equation*}
so that the first factor in \eqref{gu2} is in $\Lambda_{\ZZ_{(p)}}\lb t \rb$; the rest being as in \eqref{gup}.
\end{proof}

\subsection{Proof of \texorpdfstring{\cref{ifprop}}{Proposition 12}}\label{ifpropproof}
Recall that we assume that $\p_n(P) - \p_n(Q) \in n \aa$ for all $n$ with $1 \leq n \leq N$; we aim to show that $\e_n(P) - \e_n(Q) \in \aa$ for $n$ in the same range. We use the results of \cref{glamdef} to make some  reductions: by \eqref{englam}, it suffices to show that
\begin{equation*}
  \g_\lambda(P) - \g_\lambda(Q) \in \aa\quad\text{for $\lambda \vdash^{(p)} n$ if $1\leq n \leq N$};
\end{equation*}
by \eqref{gur} it suffices to prove that $\g_{(u)^r}(P) -\g_{(u)^r}(Q) \in \aa$ for all $u$ prime to $p$ and all~$r$ with $ur \leq N$. As in \eqref{gugessel}, write 
\begin{equation}G_u(P)(t) : = \sum_{r = 0}^\infty \g_{(u)^r}(P) \, t^{ur}\end{equation}
and the same for $Q$. By \cref{gesselprop} we know that $G_u(P)(t)$ and $G_u(Q)(t)$ are in $A \lb t \rb$; to prove the current proposition  it suffices to show that 
$$G_u(P)(t) - G_u(Q)(t) \in \aa\lb t \rb + (t^{N + 1})$$
under the assumption that $\p_{up^j} (P) - \p_{up^j} (Q) = p^{j} a_j$ for some $a_j \in \aa$ for every $j$ with $up^j \leq N$. Let $J$ be the maximum such $j$. 
We work modulo $t^{N + 1}$. Assume again for now that $p$ is odd, and again set $\eps_u = (-1)^{u-1}$. Then as in \eqref{gup} we have
\begin{align*}
\nonumber G_u(P)(t)- G_u(Q)(t) 
& =
\exp \Big( \sum_{j = 0}^\infty \eps_u \frac{ \p_{u p^j}(P)}{u p^j} \,t^{u p^j} \Big)
- G_u(Q)(t)\\
& \equiv 
\exp  \Big( \sum_{j = 0}^J\eps_u \frac{\p_{u p^j}(Q) + p^{j} a_{j}}{u p^j}\, t^{u p^j}  \Big)
- G_u(Q)(t) \mod{t^{N + 1}}.
\end{align*}
Since the exponential of a sum is the product of corresponding exponentials, we may rewrite the latter (the congruences being modulo $t^{N+1}$):
\begin{align} \label{expp2}
\nonumber G_u(P)(t) - G_u(Q)(t) 
&\equiv \exp  \Big( \sum_{j = 0}^J \eps_u \frac{\p_{u p^j}(Q) }{u p^j} \, t^{u p^j}  \Big)
\exp  \Big( \sum_{j = 0}^J \frac{ \eps_u a_{j} t^{up^j}}{u}  \Big) 
- G_u(Q)(t)\\
\nonumber
&\equiv  \exp  \Big( \sum_{j = 0}^\infty \eps_u \frac{\p_{u p^j}(Q) }{u p^j} \, t^{u p^j}  \Big)
\exp  \Big( \sum_{j = 0}^J \frac{ \eps_u a_{j} t^{up^j}}{u}  \Big) 
- G_u(Q)(t)\\
\nonumber
&=  G_u(Q)(t) \bigg(\exp  \Big( \sum_{j = 0}^J \frac{ \eps_u a_{j} t^{up^j}}{u}\Big) - 1 \bigg)\\
\nonumber  &= G_u(Q)(t) \bigg(\prod_{j = 0}^J \exp \Big(\frac{ \eps_u a_{j}}{u} \,t^{up^j} \Big)- 1 \bigg)\\
& = G_u(Q)(t)  \bigg(\prod_{j = 0}^J \Big( 1 + \sum_{k = 1}^\infty  \frac{\eps_u^k a_j^k\,  t^{ku p^j}} {u^k k! } \Big) - 1 \bigg).
\end{align}

By assumption, $\aa$ is a divided-power ideal (\cref{dpideal}), so that $a_j^k / k! \in \aa$ for every $k \geq 1$. Moreover~\mbox{$u^{-k} \in \ZZ_{(p)}$} since $u$ is prime to $p$. Therefore, for each $j$, the expression $\displaystyle\sum_{k = 1}^\infty  \frac{\eps_u^k a_j^k\,  t^{ku p^j}} {u^k k! }$ is in  $\aa\lb t\rb;$
and hence the same is true for all of 
$$\prod_{j = 0}^J \Big( 1 + \sum_{k = 1}^\infty  \frac{\eps_u^k a_j^k\,  t^{ku p^j}} {u^k k! } \Big) - 1.$$
Finally, since $G_u(Q)(t)$ in $A\lb t \rb$ (\cref{gesselprop}), we know that the last expression of \eqref{expp2}, and thus $G_u(P)(t) - G_u(Q)(t)$, is in $\aa\lb t \rb$ modulo $t^{N + 1}$, as required. 

For $p = 2$, use \eqref{gu2} in place of \eqref{gup}, so that the analogue of \eqref{expp2} is 
\begin{align*}
G_u(P)(t) - G_u(Q)(t) 
&\equiv G_u(Q)(t) \bigg(\exp  \Big( \frac{2t^u a_0}{u}\Big)  \exp  \Big(\sum_{j = 0}^J \frac{ - a_{j} t^{up^j}}{u}\Big) - 1 \bigg),
\end{align*}
again modulo $t^{N+1}$.
But the additional term $\exp  \big( \frac{2t^u a_0}{u}\big)$ is in $1 + \aa\lb t \rb$ for the same reason as $ \exp  \big(\sum_{j = 0}^J \frac{ - a_{j} t^{up^j}}{u}\big)$. 

Therefore \cref{ifprop} is proved for all primes~$p$.  \qed

\section{{Representation theory corollaries}}\label{secondproofsec}

In the case where $A$, in addition to being a torsion-free $\ZZ_{(p)}$-algebra, is a domain and the divided-power ideal $\aa$ is maximal, we can interpret a monic polynomial in $A[T]$ as the characteristic polynomial
for the action of a linear operator $T$ on a free $A$-module and the $n^{\rm th}$ power sum of its roots as the trace of $T^n$ on that module. \cref{fullalgebrathm} then becomes a statement about congruences between traces of $T^n$ implying isomorphisms between semisimplified $(A/\aa)[T]$-modules. 

We focus on the case where $A = \OO$ is a $p$-adic DVR and $\aa = \mm$ is its maximal ideal to state the following representation-theoretic version of \cref{fullalgebrathm}; \cref{theorema} is a special case. 

\begin{theorem}
\label{moduletheorem}
Let $\OO$ be a $p$-adic DVR with maximal ideal $\mm$ of ramification degree $e \leq p -1 $ and residue field $\FF$. If $M$ and $N$ are $\OO[T]$-modules, finite and free of the same rank $d$ as $\OO$-modules, then $(M \otimes \FF)^\ss \simeq (N \otimes \FF)^\ss$ as $\FF[T]$-modules if and only if for all $n$ with $1 \leq n \leq d$ we have
\begin{equation}\label{tiredofnaming} 
\tr (T^n | M) \equiv \tr (T^n | N) \pmod{n \mm}.
\end{equation}
\end{theorem}

\begin{proof}
  Let $P$ (respectively, $Q$) in $\OO[T]$ be the characteristic polynomial of the action of $T$ on $M$ (respectively, on $N$). Let $\alpha_1, \ldots,\alpha_d$ (respectively, $\beta_1, \ldots, \beta_d$) be the roots of $P$ (respectively, $Q$) in some $p$-adic DVR $\OO'$ extending $\OO$. With this notation, as detailed in~\cref{fullalgthmext}, \cref{fullalgebrathm} under the assumption $\deg P = \deg Q$ tells us that $\barr P = \barr Q$ in $\FF[X]$ if and only if $\p_n(P) \equiv \p_n(Q) \mod{n \mm}$ for all $1 \leq n \leq d$. The latter condition is equivalent to \eqref{tiredofnaming}, since $\tr (T^n | M) = \alpha_1^n + \cdots + \alpha_d^n = \p_n(P)$, and similarly $\tr(T^n|N) = \p_n(Q)$. The former condition $\barr P = \barr Q$ is equivalent to $\barr P$ and $\barr Q$ having the same multiset of roots with multiplicity in some extension of $\FF$. But the roots of $\barr P$ (respectively, $\barr Q$) are the reductions $\barr \alpha_1, \ldots, \barr \alpha_d$ (respectively, $\barr \beta_1, \ldots, \barr \beta_d$) modulo the maximal ideal $\mm'$ of $\OO'$ of $\alpha_1, \ldots, \alpha_d$ (respectively, $\beta_1, \ldots, \beta_d$). In other words, \eqref{tiredofnaming} is equivalent to the statement that, up to reordering, we have equalities
\begin{equation}\label{cactus} 
\barr \alpha_1 = \barr \beta_1,\ \barr \alpha_2 = \barr \beta_2,\ \ldots,\ \barr \alpha_d = \barr \beta_d.
\end{equation}
But the $\barr \alpha_i$ (respectively, $\barr \beta_j$) are the eigenvalues of $T$ acting on $M \otimes \FF$ (respectively $N \otimes \FF$), so that the matching in \eqref{cactus} is exactly equivalent to the up-to-semisimplification isomorphism $(M \otimes \FF)^\ss \equiv (N \otimes \FF)^\ss$. 
\end{proof}

We now return to the modular form motivation described in the introduction and prove \cref{zpquotisom}. Recall that for a $\ZZ_p$-module $M$ we write $\barr M$ for $M \otimes \FF_p$. 
\begin{corollary}[Restatement of \cref{zpquotisom}]
 \label{zpquotisom1}
  Let $M_1,M_2,N_1,N_2$ be free $\ZZ_p$-modules of finite rank, each with an action of an operator $T$.
  Suppose we have fixed $T$-equivariant embeddings $\iota_1: \barr{N_1}\into\barr{M_1}$ and $
  \iota_2: \barr{N_2}\into\barr{M_2}$ and consider the quotients
  \begin{equation*}
    W_1:=\barr{M_1}/\iota_1(\barr{N_1}),\qquad W_2:=\barr{M_2}/\iota_2(\barr{N_2}).
  \end{equation*}
  Then $W_1^\ss \simeq W_2^\ss$ as $\FF_p[T]$-modules if and only if for every $n\geq 0$ we have
  \begin{equation}\label{modncongtrace} v_p\big(\tr (T^n | M_1) -  \tr (T^n | N_1) - \tr (T^n | M_2) + \tr(T^n | N_2)\big) \geq 1 + v_p(n).
\end{equation}
\end{corollary}

\begin{proof}
Using \cref{moduletheorem}, the condition in \eqref{modncongtrace} is equivalent to the $\FF_p[T]$-module isomorphism 
\begin{equation}
\label{m1} \left(\overline{ M_1 \oplus N_2}\right)^{\ss} \simeq \left(\overline{ M_2 \oplus N_1}\right)^{\ss}.
\end{equation} Taking a quotient on the left-hand side by  $\iota_1(\barr{N_1})^\ss \oplus \barr{N_2}^\ss$ and on the right-hand side by \mbox{$\iota_2({\barr{N_2}})^\ss \oplus \barr{N_1}^\ss$} shows that \eqref{m1} is equivalent to the isomorphism $W_1^\ss \simeq W_2^\ss$. 
\end{proof} 

\begin{rem}\begin{itemize}
\item The congruence for \mbox{$0 \leq n \leq \rank M_1 + \rank N_2$} suffices in \eqref{modncongtrace}. 
\item \cref{zpquotisom1} also holds with $\ZZ_p, \FF_p, 1 + v_p(n)$ replaced by $\OO$, $\FF$, $\frac{1}{e} + v_p(n)$, respectively, where $\OO$ is a $p$-adic DVR with residue field $\FF$ and ramification degree $e \leq p - 1$ over~$\ZZ_p$.\qedhere
\end{itemize}
\end{rem}


\section*{Appendix. Brauer-Nesbitt and linear independence of characters}\label{brauernesbittsec} We briefly review the Brauer-Nesbitt theorem and connections to linear independence of characters in the setting of this paper. 

\setcounter{section}{0}

\begin{theorem}[{Brauer-Nesbitt \cite[30.16]{CurtisReiner} or \cite[Theorem 2.4.6 ff.]{wiesegaloisrep} for convenient presentation}]\label{brauernesbitt}
Let $k$ be a field; $R$ a $k$-algebra;  $V$ a semisimple $R$-module, finite dimensional as a $k$-vector space. 
  \begin{enumerate}[leftmargin=15pt]
\item {\bf Characteristic polynomial version:} The characteristic polynomial map $$r \mapsto \charpoly(r | V) \in k[X]$$ for every $r$ in~$R$ (equivalently, in a $k$-basis of $R$) determines $V$ uniquely.
\item\label{b} {\bf Trace version:} If $\cchar k= 0$ or $\cchar k > \dim_k V$ then the trace map \mbox{$r\mapsto \tr(r | V)$} for every $r$ in $R$ (equivalently,  in a $k$-basis of $R$) determines $V$ uniquely. 
\item \label{tracecomp} {\bf Trace version complement:} If $\cchar k = p$ then the trace map $r\mapsto \tr(r | V)$ for every $r$ in $R$ (equivalently,  in a $k$-basis of $R$) determines the multiplicity modulo~$p$ of every irreducible component of $V$.
\end{enumerate}
\end{theorem}

Since elementary symmetric functions determine the power-sum symmetric functions over $\ZZ$, 
the characteristic polynomial version of Brauer-Nesbitt always implies the trace version. Conversely, if $\cchar k = 0$ or $\cchar k > \dim_k V$, then $(\dim_k V)!$ is invertible in $k$, so that the power-sum functions determine the relevant elementary symmetric functions over $k$ (\cref{enplam}), and hence the trace version of Brauer-Nesbitt is equivalent to the characteristic-polynomial version. In the critical positive characteristic case $\cchar k < \dim_k V$, the trace version both assumes and concludes less than the characteristic polynomial version; neither implies the other. But if $R = k[T]$, then $R$ is abelian, so that every irreducible $R$-module is one-dimensional over $k$. In this case, both the trace version and its complement follow from the  well-known statement about linear independence of characters. 

\begin{theorem}[Linear independence of characters (Artin). See, for example, {\cite[Theorem~VI.4.1]{lang}}]\label{linindep} Let $B$ be a monoid and $\chi_1, \ldots, \chi_d: B \to k$ multiplicative characters from $B$ to a field $k$. Then $\chi_1, \ldots, \chi_r$ are $k$-linearly independent.
\end{theorem}

\begin{proposition}\label{trace} 
\cref{linindep} implies parts \ref{b} and \ref{tracecomp} of \cref{brauernesbitt}  for $R = k[T]$.
\end{proposition}

\begin{proof}
Given two finite-dimensional $k$-vector spaces $V, W$ each with the action of a single operator~$T$, let $\alpha_1, \ldots, \alpha_d$ be the list of distinct eigenvalues appearing in either $T | V$ or $T | W$ and set $B : = \ZZ^+$ and $\chi_i(n) := \alpha_i^n$. The statement that $\tr (T^n | V) = \tr (T^n | W)$ is equivalent to
$$\sum_{i = 1}^d f_i(V) \chi_i(n) = \sum_{i = 1}^d f_i(W) \chi_i(n),$$
where $f_i(V)$ is the multiplicity of $\alpha_i$ as an eigenvalue of the action of $T$ on $V$, and the same for~$W$. Linear independence of characters, then, tells us that for all~$i$ we have
$f_i(V) = f_i(W)$ in~$k$.
This simultaneously recovers for $R = k[T]$ both the trace version of Brauer-Nesbitt and its complement.
\end{proof}

The converse --- that the trace version of Brauer-Nesbitt together with its complement implies linear independence of characters --- is also true over a prime field ($k = \QQ$ or $k = \FF_p$ for some prime $p$).

\subsection*{Acknowledgments} First and foremost we thank Ira Gessel, both for his beautiful proof of \cref{gesselprop} and for allowing us to use it here. We are also grateful to Preston Wake, who patiently and generously listened to an 
error-riddled half-baked early presentation on our motivating application and both pushed and helped us to articulate the precise conditions on the ring $A$ in \cref{fullalgebrathm}. We thank John Bergdall for helpful comments. Finally we are grateful  to the Max Planck Institute for Mathematics in Bonn, whose generous hospitality allowed us to begin collaborating in 2018 and nurtured the third-named author during the Summer~2021 pandemic reprieve.

\end{document}